\documentclass[12pt]{amsart}
\usepackage{amsmath,amscd,amssymb,amsfonts,graphics}
\newcommand{\h}{\hbox}
\newcommand{\q}{\quad}

\newcommand{\bs}{\par\bigskip}
\newcommand{\ms}{\par\medskip}
\newcommand{\sk}{\par\smallskip}
\newcommand{\bsn}{\par\bigskip\noindent}
\newcommand{\msn}{\par\medskip\noindent}
\newcommand{\skn}{\par\smallskip\noindent}
\newcommand{\ges}{\geqslant}
\newcommand{\les}{\leqslant}
\newcommand{\1}{\hskip1pt}
\newcommand{\mopl}{\hbox{$\bigoplus$}}
\newcommand{\msum}{\hbox{$\sum$}}
\newcommand{\mprod}{\hbox{$\prod$}}

\newcommand{\D}{{\mathcal D}}
\newcommand{\E}{{\mathcal E}}
\newcommand{\I}{{\mathcal I}}
\newcommand{\Hc}{{\mathcal H}}
\newcommand{\M}{{\mathcal M}}
\newcommand{\OO}{{\mathcal O}}
\newcommand{\X}{{\mathcal X}}

\newcommand{\Xb}{{}\,\overline{\!X}{}}

\newcommand{\ft}{{}\,\widetilde{\!f}{}}
\newcommand{\Xt}{{}\,\widetilde{\!X}{}}
\newcommand{\Yt}{\widetilde{Y}}
\newcommand{\yt}{\widetilde{y}}
\newcommand{\DD}{{\mathbb D}}
\newcommand{\N}{{\mathbb N}}
\newcommand{\PP}{{\mathbb P}}
\newcommand{\R}{{\mathbf R}}
\newcommand{\Q}{{\mathbb Q}}
\newcommand{\C}{{\mathbb C}}
\newcommand{\Z}{{\mathbb Z}}
\newcommand{\Gr}{{\rm Gr}}
\newcommand{\IH}{{\rm IH}}
\newcommand{\al}{\alpha}
\newcommand{\alt}{\widetilde{\alpha}}
\newcommand{\Om}{\Omega}
\newcommand{\om}{\omega}
\newcommand{\De}{\Delta}
\newcommand{\bl}{\bigl}
\newcommand{\br}{\bigr}
\newcommand{\ssb}{\raise.15ex\h{${\scriptscriptstyle\bullet}$}}
\newcommand{\ssc}{\,\raise.15ex\hbox{${\scriptstyle\circ}$}\,}
\newcommand{\onto}{\twoheadrightarrow}
\newcommand{\into}{\hookrightarrow}
\newcommand{\simto}{\,\,\rlap{\hskip1.3mm\raise1.4mm\hbox{$\sim$}}\hbox{$\longrightarrow$}\,\,}
\begin{document}
\title[Rational singularities and Hodge structure]
{Deformation of rational singularities\\and Hodge structure}
\author[M. Kerr]{Matt Kerr}
\address{M. Kerr : Washington University in St. Louis, Department of Mathematics and Statistics, St. Louis, MO 63130-4899}
\email{matkerr@math.wustl.edu}
\author[R. Laza]{Radu Laza}
\address{R. Laza : Stony Brook University, Department of Mathematics, Stony Brook, NY 11794-3651}
\email{rlaza@math.stonybrook.edu}
\author[M. Saito]{Morihiko Saito}
\address{M. Saito : RIMS Kyoto University, Kyoto 606-8502 Japan}
\email{msaito@kurims.kyoto-u.ac.jp}
\thanks{During its writing, MK was supported by NSF Grant DMS-2101482 and Simons Foundation Collaboration Grant 634268; RL was supported by NSF Grant DMS-2101640; MS was supported by JSPS Kakenhi 15K04816.}
\begin{abstract} For a one-parameter degeneration of reduced compact complex analytic spaces of dimension $n$, we prove the invariance of the frontier Hodge numbers $h^{p,q}$ (that is, with $pq(n{-}p)(n{-}q)=0$) for the intersection cohomology of the fibers and also for the cohomology of their desingularizations, assuming that the central fiber is reduced, projective, and has only rational singularities. This can be shown to be equivalent to the invariance of the dimension of the cohomology of structure sheaf (which is known in the algebraizable case), since we can prove the Hodge symmetry for all the Hodge numbers $h^{p,q}$ together with $E_1$-degeneration of the Hodge-to-de Rham spectral sequence for nearby fibers, assuming only the projectivity of the central fiber.

For the proof of the main theorem, we calculate the graded pieces of the induced $V$-filtration for the first non-zero member of the Hodge filtration on the intersection complex Hodge module of the total space, which coincides with the direct image of the dualizing sheaf of a desingularization (related to Koll\'ar's conjecture on the direct images of dualizing sheaves of smooth varieties). This calculation implies also that the order of nilpotence of the local monodromy is smaller than the general singularity case by 2 in the situation of the main theorem assuming further smoothness of general fibers. We can prove a partial converse of the main theorem under some hypothesis.
\end{abstract}
\maketitle
\centerline{\bf Introduction}
\bsn
The notion of {\it rational singularity\1} was introduced by M.~Artin \cite{Ar1}, and has been studied, for instance, in the relation to {\it simultaneous resolutions\1} of versal deformations of rational surface singularities, see \cite{At,Br1,Br2,Br3,Ar2,P,Wah} among others.
The existence of a simultaneous resolution after a finite base change implies that the local monodromy is {\it semisimple\1} in the surface rational singularity case. This kind of property also appears in degenerations of hyper-K\"ahler manifolds related to smooth fillings (similar to a simultaneous resolution of surface rational singularity) generalizing Kulikov's theorem for $K3$ surfaces, see \cite{KLSV}. The above semisimplicity may be viewed as a typical case of Corollary~1 below asserting that the order of nilpotence of local monodromies is smaller than the general case by 2 (that is, at most $\dim X_t-2$) if the singular fiber $X_0$ has only rational singularities. This assertion, however, does not hold for {\it Du~Bois\1} singularities, see (2.7) below. This is the reason for which the hypothesis on the existence of {\it non-uniruled\1} component is required in \cite[Theorem~1.6]{KLSV}, whose proof uses \cite{St3} designed for Du~Bois, see also \cite{KL1}.
\sk
In examples of one-parameter degenerations of complex manifolds, we often get {\it smooth\1} total spaces (see for instance (2.7) and Remark~(2.4) below) so that the central fibers have only hypersurface singularities.
In the case of a weighted homogenous isolated hypersurface singularity with weights $(w_0,\dots,w_n)$, rational singularity is characterized by the condition that the {\it minimal exponent\1} defined by the sum of weights $\msum_{i=0}^n\,w_i$ in this case is greater than 1, see \cite[Theorem~1.11]{Wat}. This is extended to the general isolated hypersurface singularity case with minimal exponent defined as the {\it minimal spectral number\1} in \cite{St2} (see \cite{geo}), and then to the general hypersurface singularity case with minimal exponent up to sign the maximal root of the {\it reduced Bernstein-Sato polynomial\1} $b_f(s)/(s{+}1)$, see \cite[Theorem 0.4]{rat}. (The last definition of minimal exponent is compatible with the previous ones, see Remark~(1.4)\,(iv) below.)
Combining these with the Thom-Sebastiani type theorem for Bernstein-Sato polynomials \cite[Theorem 0.8]{mic}, we may have rational hypersurface singularities rather easily in the higher dimensional case.
\sk
On the other hand, Du~Bois singularity is characterized in the hypersurface case by the condition that the minimal exponent is at least 1, see \cite[Theorem~0.5]{mos}. It is well-known that this condition is equivalent to that the {\it log canonical threshold,} that is, the {\it minimal jumping coefficient\1} of the associated {\it multiplier ideals,} is 1. Note, however, that rational singularity {\it cannot\1} be characterized by using multiplier ideals, and {\it Hodge ideals\1} \cite{MP1,MP2}, or the minimal exponent as above, must be employed. It seems also possible to use Steenbrink spectrum, although we would have to calculate it at {\it every\1} point of the hypersurface near a given point {\it even locally,} see Remark~(1.4)\,(v) below.
\sk
Rational singularities are also related to birational geometry and to compactifications of moduli spaces. It is known that {\it canonical,} or more generally, {\it log-terminal\1} singularities are rational, see for instance \cite{E2,Fu,KM}, etc.
In this paper we prove the following.
\msn
{\bf Theorem~1.} {\it Let $f:X\to\De$ be a proper surjective morphism of a connected reduced complex analytic space onto a disk with relative dimension $n$. Assume the central fiber $X_0$ is reduced, projective, and has only {\it rational singularities.} Let $\Xt_t$ be a desingularization of the fiber $X_t$ $(t\in\De)$. Then the Hodge filtration $F$ on the cohomology $H^j(\Xt_t)$ and on the intersection cohomology $\IH^j(X_t)$ is strict for $t\in\De$, replacing the desingularization $\Xt_t$ and shrinking $\De$ if necessary. We have moreover the equalities
$$h_{\IH}^{p,q}(X_0)=h^{p,q}(\Xt_0)=h_{\IH}^{p,q}(X_t)=h^{p,q}(\Xt_t),
\leqno(1)$$
\vskip-7mm
$$h^q(\OO_{X_0})=h^q(\OO_{\Xt_0})=h^q(\OO_{X_t})=h^q(\OO_{\Xt_t}),
\leqno(2)$$
for $\,\,t\in\De^*$, $\,(p,q)\in\N^2\cap[0,n]^2\,$ with $\,pq(n{-}p)(n{-}q)=0\,\,\,($where $\,p=0\,$ in $\,(2))$. Here we set $\,h^{p,q}_{\IH}(X_t)\,{:=}\,\dim\Gr_F^p\IH^{p+q}(X_t)$, $\,h^{p,q}(\Xt_t)\,{:=}\,\dim\Gr_F^pH^{p+q}(\Xt_t)$, $\,h^q(\OO_{X_t})\,{:=}\,\dim H^q(X_t,\OO_{X_t})$.}
\ms
In the proof we have to assume that the union of $\Xt_t\,\,(t\,{\in}\,\De^*$) can be extended to a desingularization of the total space $X$ (forgetting about $\Xt_0$). This does not seem automatic although it is not easy to construct a counterexample.
\sk
One of our motivations for Theorem~1 comes from the study of compactifications of moduli spaces for smooth projective varieties with non-zero frontier Hodge numbers, and especially the comparison between geometric approaches (such as in \cite{KSB, Al}) and analytic approaches (based on period maps) to the construction of moduli spaces, see \cite{KL1,KL2} for more details. The assertion (2) is a refinement of the invariance of arithmetic genus under a flat deformation (\cite{Ig}), see also Remark~3 below for a related topic. We can prove a partial converse of Theorem~1 under a rather strong hypothesis, see (2.8) below. It is not yet clear whether the projectivity assumption on the central fiber can be weakened to Moishezon, see Remark~(2.5) below.
\msn
{\bf Remark 1.} It is easy to see that the {\it frontier\1} Hodge numbers $h^{p,q}(\Xt_t)$ (that is, satisfying $\,pq(n{-}p)(n{-}q)=0$) are independent of a choice of the desingularization $\Xt_t$ for any $t\in\De$ (using the Hartogs theorem (via Cauchy formula), Serre duality, and (1.1.1)), and coincide with those of the intersection cohomology $\IH^j(X_t)$ {\it when\1} $t=0$, see Corollary~(2.1) below. However, the last coincidence {\it for\1} $t\ne 0$ does not seem trivial unless $X_t$ is algebraic or $p=0,n$ (since the Hodge symmetry may fail).
The first and last equalities in (1) can be lifted to canonical isomorphisms of the corresponding vector spaces by choosing a relative ample line bundle of the desingularization morphism (using \cite{De1}, \cite{De3}, see also \cite[2.4]{toh}).
The strictness of the Hodge filtration on $H^j(\Xt_t)$ means that it is induced from the Hodge filtration on the de Rham complex which is {\it strict\1} after taking the derived global section functor. This strictness is equivalent to the $E_1$-degeneration of the Hodge-to-de Rham spectral sequence, see \cite[Proposition 1.3.2]{TH2}. For the intersection cohomology ${\rm IH}^j(X_t)$, we use the de Rham complex of the filtered $\D$-module corresponding to the intersection complex.
\msn
{\bf Remark 2.} No assumption is made for general fibers $X_t$ in Theorem~1 except for the properness of $f$ and the reducedness of $X$.
(This is changed from the first version following a suggestion of a referee.)
However, the projectivity of the central fiber $X_0$ can imply that of $X_t$ if one assumes further that $X_0$ is a variety {\it of general type,} and has only {\it canonical\1} singularities (hence only rational singularities \cite{E2}), see \cite{Kol}. Here the assumption on ``general type" is essential, since one may have a non-projective K\"ahler deformation of a smooth projective variety in the Abelian variety or $K3$ surface case.
\msn
{\bf Remark 3.} In the {\it algebraizable\1} case (that is, when $f$ is extendable to a morphism of complex algebraic varieties), the equalities in (2) can be shown by combining \cite[Lemma~1]{DJ}, \cite[Theorem 4.6]{DB} (see also \cite[Theorem 7.8]{KK}) with \cite[Theorem 2]{E1} and using the assertion that rational singularities are Du~Bois (see \cite[Corollary 2.6]{Kov}, \cite[Theorem 5.4]{mhc}), as is remarked by one of the referees.
These equalities imply (1) in the algebraizable case using Hodge symmetry.
\msn
{\bf Remark 4.} It seems quite desirable not to assume the {\it algebraizability\1} of a one-parameter degeneration as in Theorem~1. Indeed, certain analytic coordinates defined naturally near the boundary of a moduli space are not necessarily globally algebraic (for instance, in the toroidal compactification case), hence it is not necessarily easy to see whether a particular {\it analytic\1} curve given locally near the boundary can be extended globally to an {\it algebraic \1}curve on the moduli space (even in the projective one-parameter family case).
\sk
For the proof of Theorem~1, we prove the following ``Hodge symmetry", assuming only the projectivity of the central fiber (no assumption on singularities).
\msn
{\bf Proposition~1.} {\it Let $f:X\to\De$ be a proper surjective morphism of a connected reduced complex analytic space onto a disk such that the central fiber $X_0$ is projective. Then, shrinking $\De$ if necessary, the Hodge filtration $F$ on the intersection cohomology $\IH^j(X_t)$ is strict, and there are equalities $($that is, the Hodge symmetry$)$
$$h_{\IH}^{p,q}(X_t)=h_{\IH}^{q,p}(X_t)=h_{\IH}^{n-p,n-q}(X_t)\q(p,q\in\Z),
\leqno(3)$$
for any $t\in\De$. Moreover we have for $t\in\De^*$, $q\in\Z$}
$$h_{\IH}^{n,q}(X_t)=\dim H^q(X_0,F_{-n}\bl(^p\psi_f{\rm IC}_X\Q_h)\br).
\leqno(4)$$
\ms
Here ${}^p\psi_f:=\psi_f[-1]$, ${}^p\varphi_f:=\varphi_f[-1]$ are shifted nearby or vanishing cycle functors (these preserve mixed Hodge modules), and ${\rm IC}_X\Q_h$ denotes the pure Hodge module of weight $n{+}1$ whose underling $\Q$-complex is the intersection complex ${\rm IC}_X\Q$, see \cite{BBD}. For a mixed Hodge module $\M$ in general, the first non-zero member of the Hodge filtration is denoted by $F_{p(\M)}(\M)$, see (2.1) below. Note that $F_p(\M)$ is well-defined also for $p\les p(\M)$, and we have $-n\les p(\M)$ if $\M={}^p\psi_f\1{\rm IC}_X\Q_h$ or ${}^p\varphi_f\1{\rm IC}_X\Q_h$, where the strict inequality occurs in the last case, see (6) below.
\sk
Proposition~1 implies that the central fiber of a one-parameter family cannot be projective if the Hodge symmetry of nearby fibers fails. Here no assertion is claimed about the relation between $h_{\IH}^{p,q}(X_0)$ and $h_{\IH}^{p,q}(X_t)$ ($t\,{\in}\,\De^*$). Indeed, the nearby cycle complex of the intersection complex $^p\psi_f{\rm IC}_X\Q_h$ in (4) is used also for the proof of (3), and the relation between this nearby cycle complex and the intersection complex ${\rm IC}_{X_0}\Q_h$ is quite unclear (except for $F_{-n}$ under the assumption of Theorem~1 or 2 below). Proposition~1 does not necessarily hold with projectivity of $X_0$ replaced by Moishezon, see Remark~(2.4) below. (Here {\it rational singularity\1} is {\it not\1} assumed).
\sk
Applying Proposition~1 to the original $f$ and also to $\ft:\X\buildrel{\pi}\over\to X\buildrel{f}\over\to\De$ (where $\pi$ is a desingularization), the proof of Theorem~1 is reduced to the following (which is used to prove (1) only for $p=n$).
\msn
{\bf Theorem~2.} {\it Let $f$ be a holomorphic function on a reduced irreducible complex analytic space $X$. Let $Y\subset X$ be the closed analytic subspace defined by the ideal $(f)\subset\OO_X$. Assume $Y$ is reduced, and has only rational singularities. Set $n:=\dim Y$. Then $X$ has only rational singularities, replacing $X$ with an open neighborhood of $Y$, and there are isomorphisms}
$$F_{-n}({}^p\psi_f\1{\rm IC}_X\Q_h)=\om_Y,
\leqno(5)$$
\vskip-7mm
$$F_{-n}({}^p\varphi_f\1{\rm IC}_X\Q_h)=0.\,\,\,\,\,\h{}
\leqno(6)$$
\ms
We can show that the assertion (5) follows from (6) using Propositions~(1.2) and (2.1), see Remark~(2.2) below. Theorem~2 then follows from \cite[Theorem~0.6]{rat} in the $X$ smooth case.
\sk
As a corollary of Theorem~2, we can deduce the following.
\msn
{\bf Theorem~3.} {\it With the notation and assumptions of Theorem~$1$, assume $X\setminus Y$ is smooth. Set $h^{j,p,q}_{\rm lim}(X_t):=\dim\Gr_F^p\Gr^W_{p+q}H_{\rm lim}^j(X_t)$, $h^{j,p,q}_{\rm lim}(X_t)_{\ne 1}:=\dim\Gr_F^p\Gr^W_{p+q}H_{\rm lim}^j(X_t)_{\ne 1}$.
Then}
$$\aligned&h^{j,p,q}_{\rm lim}(X_t)=0\,\,\,\,\h{\it unless $(p,q)$ belongs to}\\&\q\,[1,j{-}1]^2\sqcup\{(j,0),(0,j)\}\q\q\q\q\q\q\,\,\,\,(j\les n),\\&\q\,[j{-}n{+}1,n{-}1]^2\sqcup\{(j{-}n,n),(n,j{-}n)\}\q(j>n),\endaligned
\leqno(7)$$
\vskip-3mm
$$\aligned&h^{j,p,q}_{\rm lim}(X_t)_{\ne 1}=0\,\,\,\,\h{\it if $(p,q)$ belongs to}\q\q\h{}\\&\q\{(j,0),(0,j)\}\,\,\,\,\h{\it or}\,\,\,\,\{(j{-}n,n),(n,j{-}n)\}.\q\q\q\,\,\,\h{}\endaligned
\leqno(8)$$
\ms
Here $H_{\rm lim}^j(X_t)_{\ne 1}$ is the non-unipotent monodromy part of $H_{\rm lim}^j(X_t)$ for the monodromy $T$.
Theorem~3 is a refinement of \cite[Corollary~9.9\,(ii)]{KL1}.
Since $N:=\log T_u$ is a morphism of type $(-1,-1)$ with $T=T_uT_s$ the Jordan decomposition, Theorem~3 implies the following.
\msn
{\bf Corollary~1.} {\it With the notation and assumptions of Theorem~$3$, we have a bound for the order of nilpotence of $N$ as follows\,$:$}
$$N^k=0\q\h{on}\,\,\,\,H_{\rm lim}^j(X_t)\,\,\,\,\h{for}\,\,\,\,k:=\max\bl(1,\min(j{-}1,2n{-}j{-}1)\br).
\leqno(9)$$
\ms
This bound is better than the general case by 2. (Recall that the monodromy theorem says that (9) holds in general with $\min(j{-}1,2n{-}j{-}1)$ replaced by $\min(j{+}1,2n{-}j{+}1)$. This follows from the theory of limit mixed Hodge structures, see \cite{St0}, etc.)
The assertion (9) does {\it not\1} hold for Du~Bois singularities, see (2.7) below. Corollary~1 implies that, if there is a one-parameter degeneration of smooth projective varieties such that the order of nilpotence of the local monodromy is not smaller than the upper bound in the general case by 2, then the central fiber $X_0$ cannot have only rational singularities even if we replace $X_0$ in any way.
\sk
In the case $X$ is {\it smooth\1} and $X_0$ has only rational or more generally Du~Bois singularities, we can show some relation with the {\it cohomology\1} of the singular fiber $H^j(X_0)$ which is closely related to \cite[Theorems 9.3 and 9.11]{KL1}, see Theorem~(2.6) below.
\sk
In Section 1 we review certain basics of rational singularity and its deformation (including simultaneous resolution).
In Section 2 we prove the main theorems and Theorem~(2.6) below.
\msn
{\bf Acknowledgement.} The third author was quite inspired by the arguments of \cite{KL1}, which resulted in this joint paper. We are grateful to the referees for many useful comments and questions, which led us to Proposition~1 and the current stronger version of Theorem~1 among many others.
\bs\bs
\vbox{\centerline{\bf 1. Rational singularity and smoothing}
\bsn
In this section we review certain basics of rational singularity and its deformation (including simultaneous resolution).}
\msn
{\bf 1.1.~Rational singularity.} Let $X$ be an equidimensional reduced complex analytic space. We say that $X$ has only {\it rational singularities\1} if for a desingularization $\rho:\Xt\to X$, we have the canonical isomorphism
$$\OO_X\simto\R\rho_*\OO_{\Xt},
\leqno(1.1.1)$$
or equivalently, $X$ is Cohen-Macaulay with the canonical isomorphism
$$\rho_*\om_{\Xt}\simto\om_X,
\leqno(1.1.2)$$
using duality together with the Grauert-Riemenschneider vanishing theorem \cite{GR}.
This is independent of a choice of a desingularization, since (1.1.1--2) holds for any proper morphism of complex manifolds $\rho:\Xt\to X$ inducing an isomorphism over a dense open subset (using the Hartogs theorem for (1.1.2)).
\msn
{\bf Remarks~1.1.} (i) If a reduced complex analytic space $X$ has only rational singularities, it is well-known that $X$ is {\it normal\1} and {\it Cohen-Macaulay}. Indeed, the assertion is local, and we may assume that $X$ is a closed analytic subspace of a smooth space $V$. For a desingularization $\rho:\Xt\to X$, we have ring isomorphisms $\OO_{X,x}\simto(\rho_*\OO_{\Xt})_x$ ($x\in X$) by (1.1.1) with $(\rho_*\OO_{\Xt})_x$ normal. So $X$ is normal. We may then assume $X$ globally irreducible.
\sk
Set $d_X=\dim X$. By duality for projective morphisms of complex analytic spaces together with the Grauert-Riemenschneider vanishing theorem (see Remark~(ii) below), we get the canonical isomorphisms
$$\rho_*\om_{\Xt}=\R\rho_*\om_{\Xt}=\R\rho_*(\DD\OO_{\Xt})[-d_X]
=(\DD\OO_X)[-d_X].
\leqno(1.1.3)$$
Here $\DD$ denotes the dual functor in $D^b_{\!\rm coh}(\OO_X)$, which can be defined by
$$\DD M^{\ssb}:=\tau_{\les k}{\mathcal H}om_{\OO_V}(M^{\ssb},\I^{\ssb}
)\q\bl(M^{\ssb}\in D^b_{\!\rm coh}(\OO_X)\br)\q\h{for}\q k\gg 0,
\leqno(1.1.4)$$
with $\om_V[\dim V]\simto\I^{\ssb}$ an injective resolution, if we assume that $X$ is a closed analytic subspace of a smooth space $V$ (where the dualizing sheaf $\om_V$ is canonically isomorphic to $\Om_V^{\dim V}$ using the trace morphism $H^{\dim V}_c(V,\Om_V^{\dim V})\to\C$ given by the pushdown of currents together with the Dolbeault resolution, see also \cite[2.5.1.1]{mhp}).
\sk
The equalities in (1.1.3) imply that $X$ is Cohen-Macaulay together with the canonical isomorphism (1.1.2).
\ms
(ii) Let $\rho:\Xt\to X$ be a surjective projective morphism of complex analytic spaces with $\Xt$ smooth connected and $\dim\Xt=\dim X$. Then
$$R^i\rho_*\om_{\Xt}=0\q(i>0).
\leqno(1.1.5)$$
This is known as the Grauert-Riemenschneider vanishing theorem if $\rho$ is a desingularization, see \cite[Satz 2.3]{GR} (where $X$ seems to be assumed projective, hence {\it algebraic\1}). The assertion in the {\it analytic\1} case is shown in \cite[Theorem~1]{Ta}. We can also deduce it from the {\it stability theorem\1} of polarizable Hodge modules under the direct image by a projective morphism (see \cite[Theorem~1]{mhp}) using the {\it strictness\1} of the Hodge filtration on the direct image (since the assertion is local on $X$).
\msn
{\bf 1.2.~Deformations of rational singularities.} We recall more or less well-known assertions related to deformations of rational singularities. We first see that {\it flatness\1} is satisfied in the case of one-parameter degenerations of reduced analytic spaces if and only if the total space is {\it reduced.}
\msn
{\bf Lemma~1.2.} {\it Let $f:X\to C$ be a morphism of complex analytic spaces with $\dim C=1$. Assume $C$ is smooth, $f$ is everywhere non-constant, and $X':=f^{-1}(C')$ is reduced with $C'\subset C$ the complement of a $0$-dimensional closed subspace of $C$. Then $f$ is flat if and only if $X$ is reduced.}
\msn
{\it Proof.} Since $X'$ is reduced, the kernel of the canonical surjection
$$\OO_X\onto\,\OO_{X_{\rm red}}$$
is a coherent subsheaf supported in $X\setminus X'$, and is locally annihilated by the pull-back of $t^k$ ($k\gg 0$) with $t$ a local coordinate at a point of $C\setminus C'$. So this coherent subsheaf vanishes if and only if $f$ is flat (since $\dim C=1$). This finishes the proof of Lemma~(1.2).
\msn
{\bf Proposition~1.2} (\cite{E1}). {\it Let $f$ be a non-constant holomorphic function on an irreducible reduced complex analytic space $X$. Let $Y\subset X$ be the closed analytic subspace defined by the ideal $(f)\subset\OO_X$. Assume $Y$ is reduced, and has only rational singularities. Then, replacing $X$ with a sufficiently small open neighborhood of $Y$, $X$ has only rational singularities and we have the isomorphism}
$$\om_Y=\om_X/f\om_X.
\leqno(1.2.1)$$
\msn
{\it Proof.} By definition there is a short exact sequence
$$0\to\OO_X\buildrel{f}\over\to\OO_X\to\OO_Y\to0,$$
which implies the distinguished triangle in $D^b_{\!\rm coh}(\OO_X)$
$$(\DD\OO_X)[-d_X]\buildrel{f}\over\to(\DD\OO_X)[-d_X]\to(\DD\OO_Y)[-d_Y]\buildrel{+1}\over\to.$$
Since $Y$ is Cohen-Macaulay (see Remark~(1.1)\,(i)), we have $\Hc^i(\DD\OO_Y)=0$ $(i\ne-d_Y)$. Using the associated long exact sequence together with Nakayama's lemma applied at each point of $Y$, this implies that
$$\Hc^i(\DD\OO_X)=0\q(i\ne-d_X),
\leqno(1.2.2)$$
that is, $X$ is also Cohen-Macaulay, replacing $X$ with a sufficiently small open neighborhood of $Y$. We also get the short exact sequence
$$0\to\om_X\buildrel{f}\over\to\om_X\to\om_Y\to 0,\q\h{so that}\q\om_Y=\om_X/f\om_X.
\leqno(1.2.3)$$
\sk
Let $\rho:\Xt\to X$ be a desingularization, where we may assume that $\Xt_0:=\rho^{-1}(X_0)$ is a divisor with simple normal crossings, see \cite[Ch.\,0, \S7]{Hi} and also \cite{Wlo} for a global theorem. (Note that the remaining assertion is {\it local\1} on $X$, since it is enough to show that $X$ has only rational singularities on a neighborhood of $Y$.)
More precisely, $\Xt_0\subset\Xt$ is the closed analytic subspace defined by $(\ft)\subset\OO_{\Xt}$ with $\ft:=\rho^*f$.
Let $\Yt\subset\Xt$ be the proper transform of $Y=X_0$. This is a reduced closed analytic subspace of $\Xt_0$ (of the same dimension), hence smooth.
\sk
We have the closed immersion of complex analytic spaces over $Y$
$$\Yt\into\Xt_0,$$
which induces the following morphisms using duality:
$$\rho_*\om_{\Yt}\to\rho_*\om_{\Xt_0}\to\om_Y.
\leqno(1.2.4)$$
Note that the last morphism is surjective, since the composition is.
Here we may assume that $X$, $\Xt$ are closed analytic subspaces of a polydisk $\De^m$ and $\PP^r\times\De^m$ respectively so that $\rho$ is induced by the projection $p:\PP^r\times\De^m\to\De^m$. We can use the duality isomorphism $\DD\ssc\R p_*=\R p_*\ssc\DD$ for this projection $p$. (Indeed, we have a locally free resolution whose \h{$j$\1th} component is a direct sum of copies of $\OO_{\PP^r}(-m_j)\boxtimes\OO_{{\De}^m}$ with $m_j\gg 0$ as in the algebraic case, shrinking $\De$ if necessary, see also \cite{RRV}.)
\sk
Since $\Xt\subset\PP^r\times\De^m$, we can apply the same argument as above, and get the isomorphism
$$\om_{\Xt_0}=\om_{\Xt}/\ft\om_{\Xt},$$
together with the commutative diagram
$$\begin{array}{ccccccc}
\rho_*\om_{\Xt}&\buildrel{f}\over\to&\rho_*\om_{\Xt}&\to&\rho_*\om_{\Xt_0}&\to&0\\
\downarrow&&\downarrow&&\downarrow&\raise5mm\h{}\\
\om_X&\buildrel{f}\over\to&\om_X&\to&\om_Y&\to&0
\end{array}
\leqno(1.2.5)$$
using the duality isomorphism for $\rho$ (or rather $p$) and also the Grauert-Riemenschneider vanishing theorem.
The right vertical morphism is surjective by the above argument using (1.2.4). Hence so is the middle vertical morphism by Nakayama's lemma, if we replace $X$ with a sufficiently small open neighborhood of $Y$. Proposition~(1.2) then follows.
\msn
{\bf Remarks~1.2.} (i) Proposition~(1.2) was inspired by \cite[Theorem~5.1]{Sch}, and was originally proved by assuming $X\setminus Y$ smooth. It turns out that the last hypothesis is unnecessary, and moreover the assertion in the algebraic case is already known, see \cite[Theorem~2]{E1}, where the proof is almost the same as above, and the above proof is noted simply for the convenience of the reader. A similar argument can be found also in \cite{Stv}.
\ms
(ii) If we assume $X\setminus Y$ smooth, then Proposition~(1.2) in the algebraic case is a special case of \cite[Theorem~5.1]{Sch} where the assumption that $Y$ has only rational singularities is replaced by that $Y$ has only Du~Bois singularities. (Recall that rational singularities are Du~Bois, see Remark~3 in the introduction.) It seems, however, very difficult to prove the above theorem of Schwede in an analytic setting because of the difference between rational and Du~Bois singularities.
\ms
(iii) The isomorphism of (1.2.1) {\it depends\1} on the choice of $f$ (with $Y$ fixed), although the subsheaf $f\om_X\subset\om_X$ is independent of it. (Here the division by ${\rm d}f$ is used. Note that ${\rm d}f$ trivializes the conormal sheaf of $Y\subset X$.) We have also the canonical short exact sequence
$$0\to\om_X\to\om_X(Y)\to\om_Y\to 0,
\leqno(1.2.6)$$
where the last surjection is given by residue (at least at smooth points of $Y$.) This is the dual of the short exact sequence $0\to\OO_X(-Y)\to\OO_X\to\OO_Y\to 0$ (compare with (1.2.3)).
\msn
{\bf 1.3.~Simultaneous resolution.} Let $(Y,0)$ be a germ of an isolated surface singularity. Let
$$f:X\to\De$$
be a smoothing of $Y$, that is, $f$ is flat, $X_t$ is smooth ($t\in\De^*$), and $X_0=Y$. Assume $f$ admits a {\it simultaneous resolution,} that is, there is a surjective projective morphism
$$\pi:\X\to X,$$
whose composition with $f$ is a smooth morphism $\X\to\De$ (in particular, $\X$ is smooth) and such that $\pi$ induces an isomorphism over $X\setminus X_0$.
\sk
By the commutativity of vanishing cycle functor with the direct image under a proper morphism, we have the vanishing
$$\varphi_f\R\pi_*\Q_{\X}=\R\pi_*\varphi_{\pi^*\!f}\Q_{\X}=0,
\leqno(1.3.1)$$
since $\pi^*f:\X\to\De\subset\C$ is smooth.
By the decomposition theorem for $\R\pi_*\Q_{\X}[3]$, this implies that
$$\R\pi_*\Q_{\X}[3]={\rm IC}_X\Q,
\leqno(1.3.2)$$
since the composition $\varphi_f\ssc i_*$ is the identity up to a shift of complex (using (1.3.6) below), where $i:X_0\into X$ is the inclusion. (Note that the other direct factors are contained in $X_0$). From (1.3.1--2) we thus get that
$$\varphi_f{\rm IC}_X\Q=0.
\leqno(1.3.3)$$
This is compatible with the assertion~(6) in Theorem~2. Actually we can deduce (1.3.3) from the assertion~(6) under the hypotheses of Theorem~2 assuming further that $n=2$ and the monodromy is {\it unipotent\1} (after taking a base change), since the center of symmetry of the weight filtration $W$ is shifted by 1 for the unipotent monodromy part, see (2.3.6) below. Note that $({\rm IC}_X\Q)|_{X_0}$ may change after a base change.
\sk
We can verify (1.3.3) for a concrete example as below.
\msn
{\bf Example~1.3.} Let
$$X=\{x^2+y^2+z^2=t^2\}\subset(\C^3\times\De,0),$$
where $f$ is defined by $t$, and $\De\subset\C$ is a unit disk. This is the base change of
$$g:=x^2+y^2+z^2:(\C^3,0)\to(\De,0),$$
by the cyclic double covering $(\De,0)\ni t\mapsto t^2\in(\De,0)$ (defined by using fiber product). If we blow-up $X$ at the origin, then the exceptional divisor is $\PP^1\times\PP^1$. We can blow-down this partially so that the exceptional divisor is replaced by $\PP^1$. The simultaneous resolution $\X$ can be obtained in this way, see \cite{At,Br1}, etc.
\sk
The vanishing cycle is a {\it topological\1} cycle in general, and is a sphere $S^2$ in this case. Inside the smooth family $\X=\bigsqcup_{t\in\De}\X_t$, this becomes an {\it analytic\1} cycle over $0\in\De$, and is represented by $\PP^1$. If one takes a compactification $\Xb\to\De$ of $X\to\De$ as in \cite[5.3]{SS}, then the images of the topological cycles in the compactified fiber $\Xb_t$ for $t\in\De^*$ cannot be algebraic. (Indeed, we have the surjection $F^2H^2(\Xb_t)\to H^2(X_t)$ by construction, hence the image of $H^2_c(X_t)$ in $H^2(\Xb_t)$ cannot be contained in $F^1H^2(\Xb_t)$ for $t\in\De^*$.)
So this gives an example of Hodge locus, see \cite{CDK}. In our case, however, the Hodge locus is the whole space $\De$ if we take a natural compactification, since the compactification of $\{x^2+y^2+z^2=c\}\subset\C^3$ ($c\in\De^*$) in $\PP^3$ is a smooth surface with geometric genus $p_g=0$. In this case we have a family of {\it algebraic\1} cycles over $\De$.
\sk
The normal bundle of the above $\PP^1$ in the fiber $\X_0$ at $0\in\De$ is negative. Note that the self-intersection number of the vanishing cycle is $-2$. This is related to the finiteness of the local monodromy via the {\it Picard-Lefschetz formula} \cite{Lam}, and is compatible with a criterion of analytic contraction \cite{Gr} inside $\X_0$. The contraction inside $\X$ seems more nontrivial, since the normal bundle of $\X_0\subset\X$ is trivialized by $\pi^*f:\X\to\De\subset\C$.
\sk
As for the stalk of the intersection complex ${\rm IC}_X\Q$ at $0\in X$, we have
$$\Hc^j({\rm IC}_X\Q)_0=\begin{cases}\Q&\h{if}\,\,\,j=-3,\\ \Q(-1)&\h{if}\,\,\,j=-1,\\ \,0&\h{if}\,\,\,j\ne -1,-3.\end{cases}
\leqno(1.3.4)$$
It turns out that we have the same for the nearby cycles, that is,
$$\Hc^j\psi_f({\rm IC}_X\Q)_0=\Hc^j\psi_f(\Q_X[3])_0=\begin{cases}\Q&\h{if}\,\,\,j=-3,\\ \Q(-1)&\h{if}\,\,\,j=-1,\\ \,0&\h{if}\,\,\,j\ne -1,-3,\end{cases}
\leqno(1.3.5)$$
calculating the nearby cycles $\psi_g\Q$ for $g$, since $f$ is the base change of $g$.
(Note that $g$ has an isolated singularity of type $A_1$, and has Milnor number 1, see \cite[(1.5.1)]{JKSY} for the spectrum.)
\sk
These are compatible with the vanishing of the vanishing cycles $\varphi_f{\rm IC}_X\Q$ in (1.3.3) via the long exact sequence associated with vanishing cycle triangle (see \cite{De2}):
$$i^*\to\psi_f\to\varphi_f\buildrel{+1}\over\to,
\leqno(1.3.6)$$
where $i:X_0\into X$ is the inclusion.
\sk
Note also that (1.3.4) is compatible with the short exact sequence of mixed Hodge modules
$$0\to\Q_{h,Y}[2]\to i^*({\rm IC}_X\Q_h)[-1]\to\Q_{h,\{0\}}(-1)\to 0.
\leqno(1.3.7)$$
which follows from \cite[Theorem~1]{DS1}.
\msn
{\bf Remarks~1.3.} (i) In the case of general rational surface singularities, we have a simultaneous resolution over an irreducible component (called the {\it Artin component\1}) of the base space of a miniversal deformation of a rational surface singularity after taking the base change by a ramified finite Galois covering of this component, see \cite{Ar1,Wah}, etc. Note that the covering transformation group is closely related to $-2$ curves (see \cite{Wah}, etc.), and this seems to be related to the {\it Picard-Lefschetz formula\1} as is explained below in the $A,D,E$ case.
\ms
(ii) In the rational double point case (that is, of type $A,D,E$), the base space $S$ is smooth (hence irreducible), and the covering transformation group is given by the corresponding {\it Weyl group,} see \cite{Br3}, etc.
We have a local system on the complement of the discriminant $D\subset S$ which is generated by vanishing cycles at the smooth points of the discriminant. Its monodromy group is finite, and is isomorphic to the Weyl group which is a {\it reflection group\1} generated by reflections, where a vanishing cycle at a smooth point of the discriminant determines a reflection via the Picard-Lefschetz formula as is well-known. The local system is trivialized by taking the pull-back under the unramified finite covering associated to the finite monodromy group. This triviality is needed for the simultaneous resolution, since the local system is extended over the whole base space (which is contractible) after the base change.
We can compactify the miniversal deformation $X\to S$ into a projective family over $S$ using the natural $\C^*$-action as in \cite{St1}. Here it is expected that the Hodge locus is the whole space, that is, we have a family of {\it algebraic\1} cycles over the whose base space as in Example~1.3. If we take a compactification as in \cite[5.3]{SS}, then the Hodge locus is contained in the discriminant by the same argument as in Example~1.3.
\ms
(iii) It seems interesting to examine whether the above observation in (ii) can be extended to the higher multiplicity case, where the simultaneous resolution is restricted to the Artin component. It is not very clear what happens at the other components; for instance, whether the monodromy group of the local system of vanishing cycles defined on a Zariski-open subset of an irreducible component is {\it finite\1} or not. It may be interesting to investigate this, for instance, in the case of \cite{P}.
\msn
{\bf 1.4.~Minimal exponent of hypersurfaces.} Let $X_0$ be a reduced hypersurface of a complex manifold $X$ defined by a holomorphic function $f$. We will denote $X_0$ also by $Y$. For $y\in Y$, the local {\it minimal exponent\1} $\alt_{Y,y}\in\Q_{>0}$ is defined as the maximal root of the reduced (or {\it microlocal,} see \cite{mic}) local Bernstein-Sato polynomial $b_{f,y}(s)/(s+1)$ up to sign. Globally the minimal exponent $\alt_Y$ is defined by
$$\alt_Y:=\min\bl\{\alt_{Y,y}\br\}{}_{y\in Y}.$$
Here we assume $\alt_Y$ exists by shrinking $X$ if necessary.
\sk
Let $\pi:(\Xt,\Xt_0)\to(X,X_0)$ be an embedded resolution with $E_i$ the exceptional divisors ($i\in I$) and $\Yt$ the proper transform of $Y=X_0$. We assume that $\Xt_0$ has simple normal crossings, $I$ is finite (shrinking $X$ if necessary), and $\pi$ is the composition of smooth center blow-ups (after Hironaka).
Let $m_i$, $\nu_i$ be the multiplicities of the pull-backs of $f$, $\eta$ along the exceptional divisors $E_i$ ($i\in I$) where $\eta\in\om_X$ is a local generator. Set
$$\aligned&\alt_{\pi,i}:=(\nu_i+1)/m_i,\q\alt_{\pi}:=\min\bl\{\alt_{\pi,i}\br\}{}_{i\in I},\q\alt'_{\pi}:=\min\bl\{\alt_{\pi,i}\br\}{}_{i\in I'},\\&\h{with}\q\q\q\q I':=\bl\{i\in I\,\,\big|\,\,E_i\cap\Yt\ne\emptyset\br\}\,\subset\,I.\\ \endaligned$$
\sk
The following shows some difference between rational and Du~Bois singularities in the hypersurface case (that is, the total transform is needed for Du~Bois, although the strict transform is enough for rational).
\msn
{\bf Proposition~1.4.} {\it We have the following equivalences\,$:$
\skn
{\rm(i)}\,\,\,$Y$ has at most rational singularities $\iff\alt_Y>1\iff\alt'_{\pi}>1$,
\skn
{\rm(ii)}\,\,\,$Y$ has at most Du~Bois singularities $\iff\alt_Y\ges 1\iff\alt_{\pi}\ges 1$.}
\msn
{\it Proof.} The first equivalences in (i), (ii) are shown respectively in \cite[Theorem 0.4]{rat} and \cite[Theorem 0.5]{mos}. The second equivalences are more or less well-known to specialists. We note here a short proof for the convenience of the reader.
\msn
{\it Proof of the last equivalence of {\rm (i)}.} The rationality of the singularities of $Y$ is equivalent to that
$$\nu_i-m_i>-1,\,\,\,\h{that is,}\,\,\,\alt_{\pi,i}>1\,\,(\forall\,i\in I'),\,\,\,\h{or}\,\,\,\,\alt'_{\pi}>1.
\leqno(1.4.1)$$
Indeed, the dualizing sheaf $\om_Y$ is locally generated by the ``residue" of $\eta/f$ along $Y$, which is given by using the last morphism of (1.2.6).
The singularity is rational if and only if the restriction of this residue to the smooth part of $Y$ can be extended to a holomorphic form on the proper transform of $Y$ (that is, it has no pole).
Taking the residue of the pull-back of $\eta/{\rm d}f$ along the proper transform $\Yt$ of $Y$ at $\yt\in\Yt$, we get locally
$$\eta':={\rm Res}_{\Yt}\,\pi^*(\eta/f)=u\,\mprod_{k=1}^nz_k^{\mu_k}\,{\rm d}z_1\wedge\cdots\wedge{\rm d}z_n,$$
for $u\in\OO_{\Yt,\yt}$ invertible, where $(z_0,\dots,z_n)$ is a local coordinate system of $(\Xt,\yt)$ compatible with $\Xt_0$ so that $\Yt=\{z_0=0\}$ locally, and $\mu_k:=\nu_k-m_k$ with
$$\pi^*f=v\,\mprod_{k=0}^n\,z_k^{m_k},\q\pi^*\eta=v'\,\mprod_{k=0}^n\,z_k^{\nu_k}\,{\rm d}z_0\wedge\cdots\wedge{\rm d}z_n,$$
for $v,v'\in\OO_{\Yt,\yt}$ invertible.
Note that $m_0=1$, $\nu_0=0$, since $Y$ is reduced.
The last equivalence of (i) then follows.
\msn
{\it Proof of the last equivalence of {\rm (ii)}.} Let ${\rm lct}(Y)$ be the {\it log canonical threshold\1} of a reduced hypersurface $Y$ of a complex manifold $X$. This can be defined as the minimal {\it jumping coefficient\1} of the {\it multiplier ideals\1} of $Y$, and coincides with the smallest $\al\in\Q$ such that $|f|^{-2\al}$ is not locally integrable on $X$, see \cite{Laz} (for the algebraic case). Here we shrink $X$ so that ${\rm lct}(Y)$ exists, if necessary.
The following is well-known:
$${\rm lct}(Y)=\min\{\alt_Y,1\}=\min\{\alt_{\pi},1\}\,\in\,(0,1].
\leqno(1.4.2)$$
The first equality follows for instance from \cite[Theorem 0.1]{BS}. (It is also possible to use analytic continuation in the variable $s\1$ of a functional equation associated with the Bernstein-Sato polynomial of $f$ to avoid the problem of derivation as distributions, see for instance \cite{JKSY}.)
The second equality can be verified by examining the local integrability condition for the pull-back of $f^{-\al}\eta\wedge\overline{f^{-\al}\eta}$ in terms of $\nu_i,m_i$, where $\eta\in\om_X$ is a local generator. (Recall that, for $a\in\R$, $\varepsilon\in\R_{>0}$, we have $\int_0^{\varepsilon}r^a{\rm d}r<+\infty$ if and only if $a>-1$. Here polar coordinates are used.)
\sk
We then get that
$${\rm lct}(Y)=\alt_Y=\alt_{\pi}\,\,\,\h{if one of}\,\,\,{\rm lct}(Y),\,\,\alt_Y,\,\,\alt_{\pi}\,\,\,\h{is smaller than 1}.
\leqno(1.4.3)$$
Note that the second equality does not necessarily hold without the last assumption.
This finishes the proof of the last equivalence of (ii).
\msn
{\bf Remarks~1.4.} (i) In the case of a {\it hypersurface\1} $Y$ of a complex manifold $X$, the following three conditions on the singularities of $Y$ are equivalent to each other:
\sk\q
(a) rational,\q(b) canonical,\q(c) log terminal.
\skn
The following two conditions are also equivalent:
\sk\q
(d) $Y$ : Du~Bois,\q(e) $(X,Y)$ : log canonical.
\skn
The conditions that $Y$ is respectively {\it canonical\1} and {\it log terminal\1} (which are used in birational geometry) are given by the following conditions for all $j\in J'$\,:
\sk\q
${\rm(b)}'\,\,\,\mu_j\ges0,\q{\rm(c)}'\,\,\,\mu_j>-1.$
\skn
Here we take a desingularization $\rho:\Yt\to Y$, and write
$$\om_{\Yt}\cong\rho^*\om_Y\otimes\OO_{\Yt}(\msum_{j\in J'}\,\mu_jD_j),$$
with $D_j\subset\Yt$ ($j\in J'$) the exceptional divisors and $\mu_j\in\Z$, see \cite{KS}.
In the case of a {\it pair\1} $(X,Y)$, we take an embedded resolution $\pi:(\Xt,\Xt_0)\to(X,X_0)$ with $Y=X_0$ as above so that
$$\om_{\Xt}\cong\pi^*\bl(\om_X\otimes\OO_X(Y)\br)\otimes\OO_{\Xt}(\msum_{j\in J}\,\mu_jE_j).$$
Then the condition that a pair $(X,Y)$ is {\it log canonical\1} is given by the condition
\sk\q
${\rm(e)}'\,\,\,\mu_j\ges-1\q(\forall\,j\in J).$
\skn
The above conditions are independent of the choice of a desingularization, since (logarithmic) differential forms are stable by pull-backs.
\sk
The equivalence of (a), (b), (c) follows from the definition and Proposition~(1.4)\,(i) (since $\mu_j\in\Z$). For (d), (e), we can apply Proposition~(1.4)\,(ii) or \cite[6.6]{KS} in the algebraic case.
\ms
(ii) There is a big difference between rational and Du~Bois singularities. For instance, Du~Bois singularities are not necessarily normal, and can be reducible, although they are {\it semi-normal,} see for instance \cite[Remark (i) after Proposition 5.2]{mhc} (where it is called weakly normal).
\ms
(iii) We have by \cite{MP2}
$$\alt_Y\ges\alt_{\pi}.
\leqno(1.4.4)$$
It may be possible to prove this by a microlocal version of an argument in \cite{Ka1} using an algebraic partial microlocalization as in \cite{mic}.
\ms
(iv) In the {\it isolated\1} hypersurface singularity case, the minimal exponent $\alt_Y$ coincides with the {\it minimal spectral number,} which is defined by using the mixed Hodge structure on the vanishing cohomology, see \cite{St2}. This follows by comparing \cite{Ma2} and \cite{SS,Va}.
In the {\it non-degenerate Newton boundary\1} case, the spectral numbers can be determined from the Newton polyhedron, and the minimal exponent coincides with the inverse of the minimal $c\in\Q$ such that $(c,\dots,c)$ is contained in the Newton polyhedron, see for instance \cite{exp}. In the case of a weighted homogeneous polynomial with weights $(w_0,\dots,w_n)$, the minimal spectral number coincides with the sum of weights $\msum_{i=0}^n\,w_i$, see for instance \cite[1.5]{JKSY}.
\ms
(v) In the {\it non-isolated\1} hypersurface singularity case, however, $\alt_{Y,y}$ cannot be determined by the Steenbrink spectrum {\it at\1} $y$, see \cite{rat} for the Steenbrink spectrum.
For instance, in the case of a {\it decomposable} reduced central hyperplane arrangement $Y\subset\C^4$ defined by $(x^a+y^a)(z^b+w^b)=0$ {\it with\1} $(a,b)=1$, the non-unipotent monodromy part of the vanishing cohomology at 0 vanishes.
In order to determined the local minimal exponent $\alt_{Y,y}$, we may have to calculate the Steenbrink spectrum at every $y'\ne y$ sufficiently near $y$.
\bs\bs
\vbox{\centerline{\bf 2. Proof of the main theorems}
\bsn
In this section we prove the main theorems and Theorem~(2.6) below.}
\msn
{\bf 2.1.~The first non-zero member of the Hodge filtration.} Let $\M$ be a Hodge module on a complex analytic space $X$. Let $(M_{U\into V},F)$ be the underlying filtered right $\D_V$-module associated to a local embedding $U\into V$ with $U\subset X$ open and $V$ smooth, see \cite[Remark 2.1.20]{mhp}. Note that $(M_{U\into V},F)$ is {\it canonical\1} as a consequence of \cite[3.2.6]{mhp}. (The latter implies the independence of the choice of a morphism to $V$ from a minimal local embedding of U into a smooth space). We define the {\it first non-zero member of the Hodge filtration\1} $F_{p(\M)}(\M)$ by
$$\aligned&p(\M):=\min\{p\in\Z\mid F_pM_{U\into V}\ne 0\,\,\,\h{for some}\,\,\,U\into V\},\\&F_{p(\M)}(\M)|_U:=F_{p(\M)}M_{U\into V}.\endaligned
\leqno(2.1.1)$$
These are locally independent of local embeddings $U\into V$ (since we use right $\D$-modules), and are globally well-defined. Note that $F_{p(\M)}(\M)$ is a coherent $\OO_X$-module by \cite[Lemma 3.2.6]{mhp}.
\sk
In the case of intersection complexes, we have the following.
\msn
{\bf Proposition~2.1.} {\it Let $X$ be an irreducible reduced complex analytic space of dimension $d_X$, and $\rho:\Xt\to X$ be a desingularization which is assumed to be a projective morphism. Then $p({\rm IC}_X\Q_h)=-d_X$, and we have the canonical isomorphism}
$$F_{-d_X}({\rm IC}_X\Q_h)=\rho_*\om_{\Xt}.
\leqno(2.1.2)$$
\msn
{\it Proof.} The first assertion follows from \cite[Proposition~3.2.2]{mhp} (see also \cite{kco}), since it holds on the smooth locus $X_{\rm sm}\subset X$. The last one follows from the {\it stability theorem\1} of polarizable Hodge modules under the direct images by projective morphisms (see \cite[Theorem~5.3.1 and Remark~5.3.12]{mhp}). Indeed, the latter theorem implies that the left-hand side of (2.1.2) is a direct factor of the right-hand side, using the {\it strict support decomposition\1} together with the {\it strictness\1} of the Hodge filtration on the direct image of the underlying filtered $\D$-module. So the assertion follows, since the right-hand side has no nontrivial subsheaf supported on a strictly smaller closed analytic subspace. (Here we can use also Remark~(2.1)\,(ii) below.)
This finishes the proof of Proposition~(2.1).
\msn
{\bf Corollary~2.1.} {\it In the notation and assumption of Proposition~$(2.1)$ above, assume $X$ is compact and algebraic. Then there are non-canonical isomorphisms of $\C$-vector spaces}
$$\Gr_F^p{\rm IH}^{p+q}(X)\cong\Gr_F^pH^{p+q}(\Xt)\q\h{if}\q pq(d_X{-}p)(d_X{-}q)=0.
\leqno(2.1.3)$$
\msn
{\it Proof.} Proposition~(2.1) implies the canonical isomorphisms for $p=d_X$. Corollary~(2.1) then follows using the self-duality together with the Hodge symmetry.
\msn
{\bf Remarks~2.1.} (i) In the notation and assumption of Proposition~(2.1), we have the {\it strict support decomposition\1} for pure Hodge modules
$$\Hc^j\rho_*(\Q_{h,\Xt}[d_X])=\mopl_{Z\subset X}\,\M^j_Z,
\leqno(2.1.4)$$
where $Z$ runs over irreducible closed analytic subsets of $X$, and $\M^j_Z$ is called the direct factor of $\Hc^j\rho_*(\Q_{h,\Xt}[d_X])$ with {\it strict support} $Z$ (that is, its underlying $\Q$-complex is an intersection complex supported on $Z$ with local system coefficients), see \cite[(5.1.3.5)]{mhp}. Note that $\Q_{h,\Xt}[d_X]$ denotes the pure Hodge module of weight $d_X$ whose underlying $\Q$-complex is $\Q_{\Xt}[d_X]$, and $\Hc^j\rho_*$ is defined as the direct image of a filtered $\D$-module with $\Q$-structure, see \cite[1.1]{ypg}, etc.
\ms
(ii) In the above notation and assumptions, we have by \cite[Proposition 2.6]{kco}
$$p(\M^j_Z)>-d_X\q\h{if}\q Z\ne X.
\leqno(2.1.5)$$
Note that $\M^j_Z=0$ if $Z=X$ and $j\ne 0$.
\msn
{\bf 2.2.~Proof of Theorem~2.} Set
$$\M:={\rm IC}_X\Q_h.$$
By Propositions~(1.2) and (2.1) (applied to $X$) together with \cite[(3.2.1.2), (3.2.3.1)]{mhp}, we get that
$$\aligned F_{-n-1}(\M)&=\om_X,\q\om_Y=\om_X/f\om_X\q\h{with}\\ V^{>0}\om_X&=\om_X,\q V^{>1}\om_X=f\om_X,\endaligned
\leqno(2.2.1)$$
(see also \cite{kco} for the third isomorphism). Here we use the filtration $V$ on $\om_X$ induced from the $V$-filtration of Kashiwara \cite{Ka2} and Malgrange \cite{Ma3} indexed by $\Q$ for the direct image of $\M$ by the graph embedding by $f$. This filtration is used to define the nearby and vanishing cycle functors for the underlying filtered $\D$-modules, see \cite[5.1.3.3]{mhp}.
We denote also by $V$ the induced filtration on $\om_Y=\om_X/f\om_X$. Then
$$\Gr_V^{\al}\om_Y=0\q\h{unless}\q\al\in(0,1].
\leqno(2.2.2)$$
\sk
Let ${}^p\psi_{f,1}$, ${}^p\psi_{f,\ne 1}$ be respectively the unipotent and non-unipotent monodromy part of ${}^p\psi_f$, and similarly for ${}^p\varphi_{f,1}$, ${}^p\varphi_{f,\ne 1}$.
We have the isomorphisms of mixed Hodge modules
$$\aligned{}^p\varphi_{f,1}\1\M&={\rm Coim}\bl(N:{}^p\psi_{f,1}\1\M\to{}^p\psi_{f,1}\1\M(-1)\br),\\ {}^p\varphi_{f,\ne1}\1\M&={}^p\psi_{f,\ne1}\1\M.\endaligned
\leqno(2.2.3)$$
The first isomorphism is a special case of \cite[5.1.4.2]{mhp}. Here $\,{\rm can}\,$ is surjective and $\,{\rm Var}\,$ is injective, since $\M$ has no non-trivial sub nor quotient object supported on $Y$, see for instance \cite[3.1.8]{mhp}. The second isomorphism is by the definition of the non-unipotent monodromy part of the vanishing cycle functor ${}^p\varphi_{f,\ne1}$.
\sk
The weight filtration $W$ on ${}^p\psi_f\M$, ${}^p\varphi_{f,1}\M$ are given by the monodromy filtration shifted by $n$ and $n{+}1$ respectively. We have the {\it $N$-primitive decomposition\1}:
$$\aligned\Gr^W_j{}^p\psi_f\1\M&=\mopl_{i\ges0}\,N^iP_N\Gr^W_{j+2i}{}^p\psi_f\1\M(i),\\ \Gr^W_j{}^p\varphi_{f,1}\1\M&=\mopl_{i\ges 0}\,N^iP_N\Gr^W_{j+2i}{}^p\varphi_{f,1}\1\M(i),\endaligned
\leqno(2.2.4)$$
where $P_N\Gr^W_j{}^p\psi_f\1\M$, $P_N\Gr^W_j{}^p\varphi_{f,1}\1\M$ are the $N$-primitive part defined by
$$\aligned P_N\Gr^W_{n+j}{}^p\psi_f\1\M&:={\rm Ker}\,N^{j+1}\subset\Gr^W_{n+j}{}^p\psi_f\1\M\q(j\ges 0),\\ P_N\Gr^W_{n+1+j}{}^p\varphi_{f,1}\1\M&:={\rm Ker}\,N^{j+1}\subset\Gr^W_{n+1+j}{}^p\varphi_{f,1}\1\M\q(j\ges 0),\endaligned$$
and they are 0 otherwise. Combining the first isomorphism of (2.2.3) with the primitive decompositions (2.2.4), we get that
$$P_N\Gr^W_j{}^p\psi_{f,1}\1\M=P_N\Gr^W_j{}^p\varphi_{f,1}\1\M\q\q(j\ges n+1).
\leqno(2.2.5)$$
(Note that the the weight filtration $W$ is strictly compatible with any morphism of mixed Hodge modules, and the functor assigning $\Gr^W_k$ is an exact functor of mixed Hodge modules.)
Moreover, we have by the semisimplicity of pure Hodge modules
$$\h{${\rm IC}_Y\Q_h$ is a direct factor of $P_N\Gr^W_n{}^p\psi_{f,1}\1\M\,\,(\subset{\rm Ker}\,N)$.}
\leqno(2.2.6)$$
\sk
We then get that
$$\aligned\mopl_{\al\in(0,1)}\Gr_V^{\al}\om_Y&=F_{-n}\bl({}^p\psi_{f,\ne1}\1\M\br)=F_{-n}\bl({}^p\varphi_{f,\ne1}\1\M\br),\\ \Gr_V^1\om_Y&=F_{-n}\bl({}^p\psi_{f,1}\1\M\br)\,\,\buildrel{\iota'}\over\hookleftarrow\,\,\om_Y,\endaligned
\leqno(2.2.7)$$
where the last inclusion $\iota'$ follows from (2.2.6) together with Proposition~(2.1). Indeed,
$$F_{-n}\bl(W_{n-1}{}^p\psi_{f,1}\1\M\br)=0,$$
using the $N$-primitive decomposition (2.2.4). (Note that $N$ preserves the filtration $F$ up to shift by 1, since it is defined up to sign by $\partial_tt\,{-}\,\al$ on $\Gr_V^{\al}$.)
\sk
The inclusion $\iota'$ in (2.2.7) implies the short exact sequence of coherent sheaves
$$0\to\om_Y\buildrel{\iota}\over\into\om_Y\to\E\to 0,
\leqno(2.2.8)$$
using (2.2.2) (since the filtration $V$ is decreasing). Here $\iota$ is the composition of $\iota'$ in (2.2.7) with the canonical inclusion $\iota'':\Gr_V^1\om_Y\into\om_Y$. By construction $\E={\rm Coker}\,\iota$ is a successive extension of
$$F_{-n}\bl(P_N\Gr^W_n{}^p\psi_{f,1}\1\M\br)/\iota'(\om_Y),\q F_{-n}\bl(P_N\Gr^W_k{}^p\psi_{f,1}\1\M\br)\,\,\,(k>n),$$
(using the primitive decomposition) and the direct summands of $F_{-n}\bl({}^p\psi_{f,\ne1}\1\M\br)$, that is, there is a finite filtration of $\E$ whose graded quotients are isomorphic to the above coherent sheaves. (It is enough to consider the $N$-primitive part, since $N$ preserves the filtration $F$ up to shift by 1.)
In particular, we have ${\rm Supp}\,\E\subset{\rm Sing}\,Y$. (Note that $Y_{\rm sm}\subset X_{\rm sm}$ using for instance a regular parameter system of $(Y,y)$ which can be extended to the one for $(X,y)$ by adding $f$.) These subspaces have codimension at least 2 by the normality of $Y$ so that
$$\Hc^j\DD\E\ne 0\q\h{for some}\,\,\,j\ges 2-n\,\,\,\,\h{if}\,\,\,\,\E\ne 0.$$
On the other hand, $\om_Y$ is also a Cohen-Macaulay module so that
$$\Hc^j\DD\om_Y=0\q(j\ne-n).$$
Using the long exact sequence
$$\cdots\to\Hc^{j-1}\DD\om_Y\to\Hc^j\DD\E\to\Hc^j\DD\om_Y\to\Hc^j\DD\om_Y\to\cdots,$$
associated to the dual triangle of the short exact sequence (2.2.8), we then conclude that $\E=0$, that is, $\iota$ is an isomorphism. This implies that $\iota'$, $\iota''$ are both isomorphisms (since they are injective), hence by the definitions of $\iota',\iota''$, we get that
$$F_{-n}\bl({}^p\psi_{f,1}\1\M\br)/\iota'(\om_Y)=F_{-n}\bl({}^p\psi_{f,\ne1}\1\M\br)=0.
\leqno(2.2.9)$$
So the isomorphism (5) follows. We get the vanishing of (6) from (2.2.9) using the first isomorphism of (2.2.3) (since the image of $\iota'$ is contained in the kernel of $N$, see (2.2.6)). This finishes the proof of Theorem~2.
\msn
{\bf Remark~2.2.} The assertion~(5) follows from (2.2.1) if $\Gr_V^{\al}\om_X=0$ for $\al\in(0,1)$. The last condition is satisfied if (6) holds. However, it is not easy to determine $F_{-n}({}^p\varphi_f\1{\rm IC}_X\Q_h)$ without calculating $F_{-n}({}^p\psi_f\1{\rm IC}_X\Q_h)$ in general. Note that the third isomorphism of (2.2.1) does not immediately imply the vanishing of $F_{-n}({}^p\varphi_f\1{\rm IC}_X\Q_h)$, since the filtration $F$ is {\it not\1} shifted by 1 when we define the induced filtration on $\Gr_V^0$, see \cite[(5.1.3.3)]{mhp}.
\msn
{\bf 2.3.~Proof of Theorem~3.} Set for $j\in\Z$
$$H^j_{\rm lim}:=H^j(X_0,\psi_f{\rm IC}_X\C),\q H^j_{\rm van}:=H^j(X_0,\varphi_{f,1}{\rm IC}_X\C).$$
Note that $H^j_{\rm lim}=H^j_{\rm lim}(X_t)\,\bl(:=H^j(X_0,\psi_f\C_X)\br)$, since $X\setminus X_0$ is smooth. The assertion (6) in Theorem~2 implies that
$$F^nH^j_{\rm van}=0\q\q(j\in\Z).
\leqno(2.3.1)$$
\sk
By definition the weight filtration $W$ on the nearby cycle Hodge module $^p\psi_f{\rm IC}_X\Q_h$ and the unipotent monodromy part of the vanishing cycle Hodge module $^p\varphi_{f,1}{\rm IC}_X\Q_h$ is given by the {\it {\it monodromy filtration\1} shifted by $\,n\,$ and $\,n+1\,$} respectively, see \cite[(5.1.6.2)]{mhp}.
Let $H^j_{{\rm van},1}$, $H^j_{{\rm van},\ne 1}$ be respectively the unipotent and non-unipotent monodromy part of the vanishing cohomology $H^j_{\rm van}$, and similarly for $H^j_{{\rm lim},1}$,$H^j_{{\rm lim},\ne 1}$.
The arguments in \cite[Proposition~4.2.2 and Corollary~4.2.4]{mhp} then imply that the weight filtration $W$ on $H^j_{\rm lim}$ and $H^j_{{\rm van},1}$ is given by the monodromy filtration shifted by $j$ and $j+1$ respectively. (Note that there is a serious gap at this point in certain literatures.) So there are isomorphisms
$$\aligned N^k:\Gr^W_{j+k}H^j_{\rm lim}&\simto\Gr^W_{j-k}H^j_{\rm lim}(-k)\q(k>0),\\
N^k:\Gr^W_{j+1+k}H^j_{{\rm van},1}&\simto\Gr^W_{j+1-k}H^j_{{\rm van},1}(-k)\q(k>0),\endaligned
\leqno(2.3.2)$$
The assertion for $H^j_{\rm lim}$ is compatible with the Schmid theorem (showing the coincidence of the two mixed Hodge structures). As for $H^j_{{\rm van},\ne 1}$, we have the canonical isomorphisms
$$H^j_{{\rm lim},\ne 1}=H^j_{{\rm van},\ne 1}\q(j\in\Z),
\leqno(2.3.3)$$
which follow from the canonical isomorphism $\psi_{f,\ne1}=\varphi_{f,\ne1}$.
\sk
We have the {\it $N$-primitive decomposition\1}:
$$\aligned\Gr^W_jH^j_{\rm lim}&=\mopl_{k\ges0}\,N^kP_N\Gr^W_{j+2k}H^j_{\rm lim}(k),\\ \Gr^W_jH^j_{{\rm van},1}&=\mopl_{k\ges 0}\,N^kP_N\Gr^W_{j+2k}H^j_{{\rm van},1}(k),\endaligned$$
where $P_N\Gr^W_jH^j_{\rm lim}$, $P_N\Gr^W_jH^j_{{\rm van},1}$ are the $N$-primitive part defined by
$$\aligned P_N\Gr^W_{j+k}H^j_{\rm lim}&:={\rm Ker}\,N^{k+1}\subset\Gr^W_{j+k}H^j_{\rm lim}\q(k\ges 0),\\ P_N\Gr^W_{j+1+k}H^j_{{\rm van},1}&:={\rm Ker}\,N^{k+1}\subset\Gr^W_{j+1+k}H^j_{{\rm van},1}\q(k\ges 0),\endaligned$$
and they are 0 otherwise. From the decomposition of the vanishing cycles $H^j_{{\rm lim},1}$ as in \cite[(5.1.4.2) or Corollary~4.2.4]{mhp} (together with the purity of the direct factor of $^p\Hc^{j-n}\R f_*{\rm IC}_X\Q$ supported at the origin), we can deduce the following isomorphisms and inclusions:
$$\aligned P_N\Gr^W_kH^j_{{\rm lim},1}&\,=\,P_N\Gr^W_kH^j_{{\rm van},1}\q\q(k>j+1),\\ P_N\Gr^W_{j+1}H^j_{{\rm lim},1}&\,\into\,P_N\Gr^W_{j+1}H^j_{{\rm van},1}.\endaligned
\leqno(2.3.4)$$
Note that the cokernel of the last inclusion corresponds to the direct factor of $^p\Hc^{j-n}\R f_*{\rm IC}_X\Q$ supported at the origin, see also \cite[3.1.8]{mhp}.
\sk
We have furthermore the {\it self-duality\1} isomorphisms
$$\aligned \DD H^j_{{\rm van},1}&=H^{2n-j}_{{\rm van},1}(n+1),\\ \DD H^j_{{\rm van},\ne 1}&=H^{2n-j}_{{\rm van},\ne 1}(n),\endaligned
\leqno(2.3.5)$$
where $\DD H$ denotes the dual of a mixed Hodge structure $H$. This follows for instance from \cite[(2.6.2)]{mhm} using the self-duality isomorphism
$$\DD({\rm IC}_X\Q_h)={\rm IC}_X\Q_h(n{+}1).$$
\sk
Combining (2.3.1), (2.3.5) and the Hodge symmetry, we then get
$$\aligned h_{{\rm van},1}^{j,p,q}&=0\q\h{unless}\,\,\begin{cases}p,q\in[2,j{-}1]&\h{if}\,\,\,j\les n,\\ p,q\in[j{-}n{+}2,n{-}1]&\h{if}\,\,\,j>n,\end{cases}\\
h_{{\rm van},\ne 1}^{j,p,q}&=0\q\h{unless}\,\,\begin{cases}p,q\in[1,j{-}1]&\h{if}\,\,\,j\les n,\\ p,q\in[j{-}n{+}1,n{-}1]&\h{if}\,\,\,j>n,\end{cases}\endaligned
\leqno(2.3.6)$$
with
$$h_{{\rm van},1}^{j,p,q}:=\dim_{\C}\Gr_F^p\Gr_{p+q}^WH^j_{{\rm van},1}\,\,\,\,\,(\h{similarly for}\,\,\,h_{{\rm van},\ne 1}^{j,p,q}).$$
Theorem~3 now follows from (2.3.6) using (2.3.2--4). This finishes the proof of Theorem~3.
\msn
{\bf 2.4.~Proof of Proposition~1.} Consider the nearby cycle Hodge module
$$^p\psi_f{\rm IC}_X\Q_h.$$
Its weight filtration is the {\it monodromy filtration shifted by\1} $n$. Hence its cohomology
$$H^j\bl(X_0,{}^p\psi_f{\rm IC}_X\Q_h\br)
\leqno(2.4.1)$$
is a mixed Hodge structure such that the weight filtration is the monodromy filtration shifted by $n{+}j$, see \cite[5.3.5]{mhp}. The cohomology (2.4.1) is defined by using $H^j(a_{X_0})_*$ with $a_{X_0}:X_0\to pt$ a canonical morphism, and this direct image is defined as that of a {\it filtered $\D$-module with $\Q$-structure,} see \cite[1.1]{ypg}, etc. (Note that $X_0$ is a projective variety, and is globally embeddable into a smooth variety $\PP^N$.) We have the {\it strictness\1} of the Hodge filtration $F$ of
$$(a_{X_0})_*{}^p\psi_f{\rm IC}_X\Q_h.$$
We have a similar assertion for the vanishing cycle Hodge module with unipotent monodromy $^p\varphi_{f,1}{\rm IC}_X\Q_h$. Moreover we can prove the direct sum decomposition in \cite[(5.1.4.2)]{mhp} using \cite[4.2.4]{mhp}.
\sk
We can extend \cite[3.3.17]{mhp} (and \cite[5.3.4]{mhp}) by \cite[2.3.9]{mhp} (with $U'_i=V'_i=\De$ for any $i$) extending \cite[\S3]{ana} to the filtered case. Here we have to use the canonicity of $(M_{U\into V\times\C},F)$ explained in (2.1).
We then get the {\it strictness\1} of the Hodge filtration $F$ on the underlying complexes of filtered $\D$-modules of
$$f_*\1{\rm IC}_X\Q_h,$$
shrinking $\De$ if necessary. (Here the direct image is defined as that of a filtered $\D$-module with $\Q$-structure using \cite[2.3.9]{mhp} as is explained above.)
Moreover the mixed Hodge structure (2.4.1) is identified with
$$^p\psi_t\Hc^jf_*{\rm IC}_X\Q_h.
\leqno(2.4.2)$$
This means that we have the ``limit mixed Hodge structure" (2.4.1--2) of the filtered $D$-modules with $\Q$-structure
$$(\Hc^jf_*{\rm IC}_X\Q_h)|_{\De*},
\leqno(2.4.3)$$
such that the weight filtration is the monodromy filtration shifted by $n{+}j$. Here it is totally unclear whether these are really {\it variations of Hodge structure\1} on $\De^*$ (although {\it polarized mixed Hodge structures\1} in the strong sense correspond to {\it nilpotent orbits,} see for instance \cite[3.13]{CKS}). Nevertheless they gives filtered free sheaves on $\De$ such that the graded quotients are free, by taking their intersections with the $V$-filtration $V^{>0}$ instead of the restriction to $\De^*$, see \cite[(3.2.3.1)]{mhp} and also the proof of \cite[3.1.3]{mhp}. Here freeness is equivalent to torsion-freeness, since $\dim\De=1$. (We may assume that the local monodromy is unipotent by taking an appropriate base change.) This freeness implies that the dimensions of the graded quotients of the Hodge filtration
$$\dim\Gr_F^p{\rm IH}^{p+q}(X_t)\q(t\in\De^*)
\leqno(2.4.4)$$
do not change by passing to the ``limit mixed Hodge structure" (2.4.1--2). This implies (4).
\sk
By the above argument the stalks at $t\in\De^*$ of the filtered $D$-module with $\Q$-structure in (2.4.3) is identified with
$${\rm IH}^{n+j}(X_t),
\leqno(2.4.5)$$
shrinking $\De$ if necessary. So we get the assertion on the strictness of the Hodge filtration for $t\in\De^*$. (The assertion for $t=0$ follows from \cite[5.3.1]{mhp}.)
\sk
Using the monodromical property of the weight filtration, the hard Lefschetz theorem, and the self-duality of the nearby cycles, we now get the symmetries of the Hodge numbers in (2.4.4) with respect to the two involutions $\iota_1,\iota_2$ of $\Z{\times}\Z$ defined by
$$\iota_1(p,q)=(q,p),\q\q\iota_2(p,q)=(n{-}p,n{-}q).
\leqno(2.4.6)$$
(For $\iota_1$ we can also employ \cite[3.13]{CKS} together with \cite[4.2.2]{mhp}. We need the self-duality for $\iota_2$.)
These symmetries for the intersection cohomology with $t=0$ follow from \cite[5.3.1]{mhp} using the projectivity of $X_0$. So Proposition~1 follows.
\msn
{\bf Remark~2.4.} Proposition~1 does not hold if the assumption on the projectivity of $X_0$ is replaced by Moishezon. Indeed, there are examples constructed by Kodaira \cite{Kod} (which is explained in \cite{Cl}) and Oda \cite[14.3]{O}. Here the total space $X$ is smooth, the general fibers $X_t$ are Hopf or Inoue surfaces, and the central fiber $X_0$ is obtained by identifying two disjoint smooth rational curves on a smooth rational surface $S$. The latter surface is obtained by successive one-point blow-ups of $\PP^1{\times}\PP^1$ and one blow-down in the Inoue surface case, while the blow-up is once in the Hopf surface case (since $S$ is a blow-up of $\PP^2$). The surface contains four or five non-singular rational curves $C_i$ $(i\in[1,4]$ or $[1,5]$) intersecting transversally, and the complement of their union is a torus. Their intersection matrix is given by
$$\scalebox{0.9}{$\begin{pmatrix}0&1&0&1\,\\1&-1&1&0\,\\0&1&0&1\,\\1&0&1&1\,\end{pmatrix}$}\q\h{or}\q\scalebox{0.9}{$\begin{pmatrix}-2&1&0&0&1\\1&-1&1&0&0\\0&1&-1&1&0\\0&0&1&1&1\\1&0&0&1&\,0\,\end{pmatrix}$}$$
Here $C_2$ and $C_4$ are identified in order to get $X_0$. We can verify that there is no linear combination $D=\msum_i\,a_iC_i$ ($a_i\in\Z$) such that $(D,C_i)>0$ for any $i$, and $(D,C_2)=(D,C_4)$. So $X_0$ cannot be projective, considering the pull-back to $S$ of a possible ample line bundle on $X_0$.
Note that this non-projectivity also follows from Proposition~1, since the Hodge symmetry fails for Hopf or Inoue surfaces.
\sk
For the calculation of $H^j(X_0)$ and the limit mixed Hodge structure $H^j_{\rm lim}$, let $\rho:S\to X_0$ be the canonical morphism. Put $C:={\rm Sing}\,X_0\,(\cong\PP^1)$. We have a short exact sequence
$$0\to\Q_{X_0}\to\rho_*\Q_S\to\Q_C\to0,
\leqno(2.4.7)$$
and
$$\Gr^W_k{}^p\psi_f(\Q_X[3])=\begin{cases}\Q_C[1]&(k=1),\\ \rho_*\Q_S[2]&(k=2),\\ \Q_C(-1)[1]&(k=3),\end{cases}
\leqno(2.4.8)$$
where these vanish otherwise. In the Hopf surface case, we then get that
$$H^j(X_0)=\begin{cases}\Q&(j=0,1),\\ \Q(-1)&(j=2),\\ \Q(-2)&(j=4),\end{cases}\q\q H^j_{\rm lim}=\begin{cases}\Q&(j=0,1),\\ \,0&(j=2),\\ \Q(-2)&(j=3,4),\end{cases}
\leqno(2.4.9)$$
with $H^{1+2i}_{\rm van}=\Q(-i)$ ($i=0,1)$, and these vanish otherwise.
In the Inoue surface case, the dimensions of the terms in (2.4.9) are increased by 1 for $j=2$.
This calculation is compatible with the assertion that the general fibers $X_t$ ($t\ne 0$) are Hopf or Inoue surfaces, for which Hodge symmetry does not hold (indeed, $h^{1,0}(X_t)=0$, $h^{0,1}(X_t)=1$). Note that the monodromy is unipotent, since $X$ is smooth and $X_0$ is a reduced divisor with normal crossings. So the local invariant cycle theorem does not hold for $H^3_{\rm lim}$.
\msn
{\bf 2.5.~Proof of Theorem~1.} Let $\pi:\X\to X$ be a desingularization such that $\X_0:=\pi^{-1}X_0$ is a divisor with normal crossings. Here we may assume that $\pi$ is projective. This implies that $\X_0$ is projective, since $X_0$ is. Let $\pi_t:\X_t\to X_t$ be the restriction of $\pi$ over $t\in\De$. This is a desingularization for $t\in\De^*$ (shrinking $\De$ if necessary). We then put $\Xt_t:=\X_t$ ($t\in\De^*$).
\sk
By Proposition~(1.2), we may assume that $X$ has only rational singularities shrinking $\De$ if necessary. Then the fibers $X_t$ ($t\in\De$) have also only rational singularities. This can be seen by taking the direct image of the short exact sequence
$$0\to\OO_{\X}\buildrel{\!\ft-c}\over\longrightarrow\OO_{\X}\to\OO_{\X_t}\to 0.$$
Here $\ft:=f\ssc\pi$, and $c$ denotes the value of the coordinate of $\De\subset\C$ at $t$.
\sk
We apply Proposition~1 to $f$ and also to $\ft$ (where ${\rm IC}_{\X}\Q_h=\Q_{h,\X}[n{+}1]$). This implies the Hodge symmetry of the terms in (1) for any $t\in\De$ (since it is trivial for $h^{p,q}(\Xt_0)$).
The proof of (1) is then reduced to the case $p=n$, and follows from (4), Theorem~2 and Proposition~(2.1).
The equalities in (2) then follow from (1) together with (1.1.1), since rational singularities are Du~Bois so that
$$\Gr_F^0H^j(X_t)=H^j(X_t,\OO_{X_t})\q\q(t\in\De),$$
see \cite[Corollary 2.6]{Kov}, \cite[Theorem 5.4]{mhc}. This finishes the proof of Theorem~1.
\msn
{\bf Remarks~2.5.}\,(i) It is not quite clear whether the assumption on the projectivity of $X_0$ in Theorem~1 can be replaced by Moishezon. One problem is that it is quite nontrivial whether we have a desingularization of the total space $\pi:\X\to X$ such that the central fiber $\X_0:=\pi^{-1}(X_0)$ is {\it projective.} This is false without the hypothesis on rational singularities, see Remark~(ii) just below. The situation is, however, rather unclear if we assume $X_0$ is {\it normal,} in particular, it has no ``self-intersection" as in Remark~(2.4) (here the {\it proper transform\1} of $X_0$ by the blow-up of the total space $\X$ along the singular locus of $X_0$ is {\it projective,} and the {\it non-connectedness\1} of its intersection with the exceptional divisor is the cause of problem).
\ms
(ii) If we do not assume that $X_0$ has only rational singularities, there are examples such that the {\it total transform\1} $\X_0:=\pi^{-1}(X_0)$ is {\it not\1} projective for any embedded desingularization $\pi:(\X,\X_0)\to(X,X_0)$. Indeed, it is enough to take an example such that the total space $X$ is smooth and the local invariant cycle theorem does not hold as in \cite{Cl} (see also Remark~(2.4) above). Here $\Q_X$ is a direct factor of $\R\pi_*\Q_{\X}$ using the self-duality of $\Q_X[n]$, $\Q_{\X}[n]$. It is known that the local invariant cycle theorem is equivalent to the decomposition theorem if the total space is smooth and the base space is a non-compact connected smooth curve, see for instance \cite[Theorem A.4]{KL1}. Note also that the projectivity of $X_0$ and the smoothness of $X$ imply the local invariant cycle theorem.
\ms
(iii) The constancy of $\dim H^j(X_t,\OO_{X_t})$ and $\dim\Gr_F^0H^j(X_t)$ ($t\in\De$) holds in the case of two examples in Remark~(2.4) (although it fails for $\dim\Gr_F^2H^3(X_t)$). It may be possible that this constancy is true under the hypothesis that the special fiber $X_0$ is Moishezon and Du~Bois (forgetting about the other frontier Hodge numbers, that is, $\dim\Gr_F^pH^{p+q}(X_t)$ with $q(n{-}p)(n{-}q)=0$, $p\ne 0$) as is questioned by a referee.
\msn
{\bf 2.6.~Case $X$ smooth.} Under a strong assumption that $X$ is {\it smooth,} we can deduce certain relations with the cohomology of the singular fiber. This is closely related to \cite[Theorems 9.3 and 9.11]{KL1}, see Remark~(2.6) below.
\msn
{\bf Theorem~2.6.} {\it Let $f:X\to\De$ be a surjective projective morphism of a complex manifold onto a disk such that general fibers $X_t\,\,\,(t\in\De^*)$ are smooth and connected, and the singular fiber $Y:=X_0$ is reduced. Let $\rho:\Yt\to Y$ be a desingularization.
In the notation of $(2.3)$, we have the following.
\skn
{\rm (a)\,} If $\,Y$ has only rational singularities, then
$$\aligned\Gr_F^pH^j(Y)&=\Gr_F^pH^j_{\rm lim}=\Gr_F^pH^j_{{\rm lim},1}=\Gr_F^p\Gr^W_jH^j_{{\rm lim},1}\\&=\Gr_F^pH^j(\Yt)\q\q\q\q(j\in\Z,\,\,p=0,n),\endaligned
\leqno(2.6.1)$$
$$\Gr_F^1H^j(Y)=\Gr_F^1H^j_{{\rm lim},1}\q\q(j\in\Z).
\leqno(2.6.2)$$
\skn
{\rm (b)\,} If $\,Y$ has only Du~Bois singularities, then}
$$\Gr_F^0H^j(Y)=\Gr_F^0H^j_{{\rm lim},1}=\Gr_F^0H^j_{\rm lim}\q\q(j\in\Z).
\leqno(2.6.3)$$
\msn
{\it Proof.} Since $X$ is smooth, we have a short exact sequence of mixed Hodge modules
$$0\to\Q_{h,Y}[n]\to{}^p\psi_{f,1}(\Q_{h,X}[n{+}1])\to{}^p\varphi_{f,1}(\Q_{h,X}[n{+}1])\to 0,
\leqno(2.6.4)$$
inducing the vanishing cycle sequence (see \cite{De2})
$$\to H^{j-1}_{{\rm van},1}\to H^j(Y)\to H^j_{{\rm lim},1}\to H^j_{{\rm van},1}\to.
\leqno(2.6.5)$$
\sk
We first show the assertion (a) assuming $Y$ has only rational singularities. The first isomorphism of (2.6.1) and (2.6.2) follow from (2.3.6) and (2.6.5). The other isomorphisms of (2.6.1) follow from Theorems~1 and 3. This finishes the proof of the assertion (a).
\sk
For the assertion (b), assume $Y$ has only Du~Bois singularities. Then
$$F_{-n}\bl({}^p\varphi_{f,\ne 1}(\Q_{h,X}[n{+}1])\br)=0,$$
by \cite[Theorem 0.5]{mos} or \cite[4.3]{MSS}. This implies (2.3.6)$_{\ne 1}$, that is, (2.3.6) holds for the non-unipotent monodromy part, in the Du~Bois case.
As to the unipotent monodromy part, it follows from \cite[(3.2.3.1), (5.1.3.3)]{mhp} that
$$F_{-n-1}\bl({}^p\varphi_{f,1}(\Q_{h,X}[n{+}1])\br)=0.$$
Using the isomorphisms in (2.3.2) and the self-duality (2.3.5), this implies that the following holds {\it unconditionally\1}:
$$h_{{\rm van},1}^{j,p,q}=0\q\h{unless}\,\,\begin{cases}p,q\in[1,j]&\h{if}\,\,\,j\les n,\\ p,q\in[j{-}n{+}1,n]&\h{if}\,\,\,j>n.\end{cases}
\leqno(2.6.6)$$
\sk
From (2.3.3), (2.3.6)$_{\ne 1}$, (2.6.6), we can now deduce the isomorphisms in (2.6.3). This finishes the proof of Theorem~(2.6).
\msn
{\bf Remark~2.6.} The assertion (a) improves \cite[Theorem~9.11]{KL1}, and the assertion (b) proves \cite[Theorem~9.3]{KL1} without assuming $X$ extendable to an algebraic variety, but assuming $X$ {\it smooth}.
\msn
{\bf 2.7.~An example of Du~Bois singularity.} Let $X\subset\PP^{n+1}\times\De$ ($n\ges 2$) be a flat family of projective hypersurfaces of degree $n{+}2$ over $\De$ defined by the equation
$$\msum_{i=1}^{n+1}\,x_i^{n+2}+x_1\cdots x_{n+2}=x_{n+2}^{n+2}\1t.$$
Here the total space is smooth shrinking $\De$ if necessary. (In general, a subvariety of $\PP^{n+1}\times\De$ defined by $f+gt=0$ for two homogeneous polynomials $f,g$ of the same degree is smooth (shrinking $\De$ if necessary) provided that $\{f=0\}\subset\PP^{n+1}$ has only isolated singularities which do not intersect $\{g=0\}$.)
\sk
Restricting to the affine space $\{x_{n+2}\ne 0\}=\C^{n+1}$, the singular fiber $X_0|_{\{x_{n+2}\ne 0\}}$ is defined by
$$h:=\msum_{i=1}^{n+1}\,y_i^{n+2}+y_1\cdots y_{n+1}=0,$$
where $y_i:=x_i/x_{n+2}$. This has an isolated Du~Bois singularity at 0, since the minimal spectral number is 1, see \cite{exp}, \cite[Theorem 0.5]{mos} (that is, the first equivalence in Proposition~(1.4)\,(ii)).
Note that $X_0$ is a rational variety. In the case $n=2$, this is closely related to \cite[Theorem II]{Kul} (which is partially generalized to \cite[Theorem 1.6]{KLSV}).
\sk
We have the following for $n\ges 2$\,:
$$N^n\ne 0\,\,\,\h{on}\,\,\, H^n_{{\rm lim},1},\,\,\h{that is, the order of nilpotence is n+1,}
\leqno(2.7.1)$$
where $H^n_{{\rm lim},1}$ is the unipotent monodromy part of the limit mixed Hodge structure.
This gives examples such that Corollary~1 does not hold if we replace rational with Du~Bois.
\ms
The assertion (2.7.1) is equivalent to that
$$N^{n-1}\ne 0\,\,\,\h{on}\,\,\, H^n_{{\rm var},1},
\leqno(2.7.2)$$
since $N^n={\rm Var}\ssc N^{n-1}\ssc{\rm can}$ on $H^n_{{\rm lim},1}$, where $H^n_{{\rm var},1}$ is the unipotent monodromy part of the vanishing cohomology.
This is analogous to \cite{Ma1} for the non-unipotent monodromy part,
see \cite[Example 3.16]{St2} for the case $n=2$, and \cite{Sta}, \cite[Corollary~3]{des} for $n\ges 2$.
(This seems to be a special case of a theory on Milnor monodromies of Newton non-degenerate functions, although certain arguments do not seem necessarily easy to follow in some papers quoted in \cite{Sta}. Note also that one cannot determine the order of nilpotence from the {\it combinatorial data\1} of an embedded resolution in general without assuming Newton non-degeneracy, see for instance \cite{DS2}.)
\msn
{\bf 2.8.~Partial converse of Theorem~1.} If $X\,{\setminus}\,Y$ is smooth, $Y$ has at most Du~Bois singularities, and the non-rational locus of $Y$ is discrete (for instance, if $Y$ has only isolated Du~Bois singularities), then we have a partial converse of Theorem~1 in the {\it algebraic\1} case as follows.
\msn
{\bf Proposition~2.8.} {\it In the notation of Theorem~$1$, assume $X$ can be extended to a complex projective variety, $X\,{\setminus}\,Y$ is smooth, $Y$ has at most Du~Bois singularities which are rational outside a finite number of points, and $(1)$ in Theorem~$1$ holds. Then $Y$ has at most rational singularities everywhere.}
\msn
{\it Proof.} The hypotheses imply that $X$ has only rational singularities \cite[Theorem 5.1]{Sch}.
In the notation of (2.2), we then get
$$F_{-n-1}({\rm IC}_X\Q_h)=\om_X,\q F_{-n}({\rm IC}_Y\Q_h)=\rho_*\om_{\Yt}\subset\om_Y=\om_X/f\om_X,
\leqno(2.8.1)$$
with $\rho:\Yt\to Y$ a desingularization. (Note that a locally principal divisor on a Cohen-Macaulay variety is Cohen-Macaulay, see for instance \cite[II, Theorem 8.21A]{AG}.) By the arguments in (2.2), we have the isomorphisms
$$\aligned\Gr_F^nH^{j+n}(\Yt)&=H^j(\Yt,\om_{\Yt}),\\
\Gr_F^nH^{j+n}_{\rm lim}(X_t)&=\mopl_{\al\in(0,1]}\,H^j(Y,\Gr_V^{\al}\om_Y)\q\q(j\in\Z),\endaligned
\leqno(2.8.2)$$
where the filtration $V$ on $\om_Y=\om_X/f\om_X\,(=V^{>0}\om_X/V^{>1}\om_X$) is as is explained in (2.2). So the equalities in (1) for $p=n$ imply that
$$\E':=\om_Y/\rho_*\om_{\Yt}=0,
\leqno(2.8.3)$$
using the Grauert-Riemenschneider theorem (which implies that $\R\rho_*\om_{\Yt}=\rho_*\om_{\Yt}$) together with the {\it Euler characteristic,} where we can deduce $\E'=0$ from $\chi(Y,\E')=0$, since the support of the $\OO_Y$-module $\E'$ is discrete. This finishes the proof of Proposition~(2.8).
\msn
{\bf Remarks~2.8.} (i) If we do not assume that the non-rational locus is discrete, we cannot conclude that $\E'=0$. For instance, if $\E'\cong i_*\OO_{\PP^1}(-1)$ with $i:\PP^1\into Y$ a closed immersion, we have
$$H^j(Y,\E')=0\q(\forall\,j\in\Z).$$
\ms
(ii) There may be a counterexample to the converse of Theorem~1 if we do not assume $Y$ Du~Bois. In the $Y$ Du~Bois case, this seems to be a quite nontrivial question.

\end{document}